\documentclass[reqno]{amsart}

\usepackage{amscd}
\usepackage{amssymb}
\usepackage{youngtab}

\theoremstyle{plain}
\newtheorem{theorem}{Theorem}[section]
\newtheorem{lemma}[theorem]{Lemma}
\newtheorem{proposition}[theorem]{Proposition}
\newtheorem{corollary}[theorem]{Corollary}

\theoremstyle{definition}
\newtheorem{definition}[theorem]{Definition}
\newtheorem{def-prop}[theorem]{Definition-Proposition}
\newtheorem{example}[theorem]{Example}

\theoremstyle{remark}
\newtheorem{remark}{Remark}[section]
\begin{document}
\title[Quiver varieties]{Some remarks on Nakajima's quiver varieties of type {\it A}.}
\author{D.A.Shmelkin}
\address{117437, Ostrovitianova, 9-4-187, Moscow, Russia.}
\email{mitia@mccme.ru}
\begin{abstract}
We try to clarify the relations between quiver varieties of type A and Kraft-Procesi proof of
normality of nilpotent conjugacy classes closures.
\end{abstract}
\maketitle
\section{Introduction}
Kraft and Procesi proved in \cite{kp} that for any nilpotent $n\times n$ matrix $A$ over an algebraically
closed field ${\bf k}$ of charactersistic zero, the closure $\overline{C_A}$ of the conjugacy
class $C_A$ of $A$ is normal, Cohen-Macaulay with rational singularities. The main idea of the proof of this
wonderful theorem is as follows: $\overline{C_A}$ is proved to be isomorphic to the categorical quotient for
an affine variety $Z$ of representations of a quiver with relations: $\overline{C_A}\cong Z/\!\!/H$,
where $H$ is a reductive group. Moreover, this $Z$ is proved to be a reduced irreducible normal complete
intersection, and this implies all the claimed properties of $\overline{C_A}$  as being inherited by
the categorical quotients over reductive groups in general.

Nakajima in \cite{na94} and \cite{na98} introduced a setup related to the term {\it quiver variety}.
A very particular case of that setup, when the underlying quiver is of type {\it A} and
the additional vector spaces are of special dimension vector leads to the above variety $Z$ used
by Kraft and Procesi. Nakajima employed this observation in \cite{na94} to illustrate quiver varieties, 
in particular, he proved a nice theorem (\cite[Theorem~7.3]{na94}) relating
the quiver variety in this case with the cotangent bundle over a flag variety. The proof is based
on another result (\cite[Theorem~7.2]{na94}) that he claimed to be 
proved in \cite{kp}. Actually, that result was proved in  \cite[Proposition~3.4]{kp} only
for special dimension vectors, not in the generality needed for Theorem 7.3. 
Unfortunately, this confusion haven't been corrected so far and we want to fill this gap, 
and without any contradiction with the valuable {\it sense} of Nakajima's result.

First of all, both Theorems 7.2 and 7.3 are true and we give proofs for them.
In addition, we show that Theorem 7.3 is closely related with a result on $\Delta$-filtered modules of Auslander algebra from \cite{bhrr}.
On the other hand, the main part of the results of \cite{kp}
(because \cite[Proposition~3.4]{kp} is only a small part of these) 
can not be generalized, in particular, 
the variety $Z$ can be reducible (see Example \ref{reduc}).

Our study does not claim to be a new result. Quite the contrary, we are trying to present 
the known results in their uncompromising beauty. 
\section{Kraft-Procesi setup and Nakajima's Theorem 7.2.}
We present the setup used in \cite{kp} keeping the local notation.
Consider a sequence of $t$ vector spaces and linear mappings between them:
\begin{equation}
U_1 \overset{A_1}{\underset{B_1}\rightleftarrows}  U_2 
\overset{A_2}{\underset{B_2}\rightleftarrows} U_3
\cdots
U_{t-1} \overset{A_{t-1}}{\underset{B_{t-1}}\rightleftarrows} U_t
\end{equation}
\noindent Consider moreover the equations as follows:
\begin{equation}\label{def_Z}
B_1A_1 = 0; B_2A_2 = A_1B_1; B_3A_3 = A_2B_2;\cdots;B_{t-1}A_{t-1} = A_{t-2}B_{t-2}
\end{equation}
\noindent and denote by $Z$ the closed subvariety defined by these equations.
The equations can be thought of as "commutativity" conditions for every $i=2,\cdots,t-1$:
two possible compositions of $U_{i-1} \rightleftarrows U_i \rightleftarrows U_{i+1}$
yield the same endomorphism of $U_i$. The extra condition $B_1A_1=0$ combined with that commutativity
implies $(A_1B_1)^2=A_1(B_1A_1)B_1=0$. Inductively, we have for $i=2,\cdots,t-1$:
\begin{equation}
(B_iA_i)^i=(A_{i-1}B_{i-1})^i= A_{i-1}(B_{i-1}A_{i-1})^{i-1}B_{i-1}=0
\Rightarrow (A_iB_i)^{i+1}=0
\end{equation}
\noindent so all these endomorphisms are nilpotent.
Denote $\dim U_i$ by $n_i$; so we have the dimension vector $(n_1,\cdots,n_t)$.
The variety $Z$ is naturally acted upon by the group $G= GL_{n_1}\times\cdots\times GL_{n_t}$
and its normal subgroup $H=GL_{n_1}\times\cdots\times GL_{n_{t-1}}$.
The above setup is interesting
for any dimension vector but each of the texts \cite{kp} and \cite[\S 7]{na94} considered those 
important for their purposes. Nakajima considered (in slightly different notation) {\it monotone} dimension vectors, that is, subject to the condition $n_1 < n_2 < \cdots < n_t$. One of the statements we feel necessary to
clarify is the following (in our reformulation consitent with given notation):
\begin{theorem}\label{7.2}
{\rm (Theorem 7.2 from \cite{na94})} Assume $(n_1,\cdots,n_t)$ is monotone.
Then the map $(A_1,B_1,\cdots,A_{t-1},B_{t-1})\to A_{t-1}B_{t-1}:Z\to {\rm End}(U_t)$ 
is the categorical quotient with respect to
$H$ and the image is the conjugacy class closure for a nilpotent matrix.
\end{theorem}

Instead of the proof for this Theorem it is stated in \cite{na94} that this result is proved
in \cite{kp}. This is not true, because in \cite{kp} a smaller subset of dimensions was considered
and the most part of the results concerns this subset, though the developed methods do allow to recover
the proof of the above Theorem (see \S\ref{pf}).

For the main goal of \cite{kp} it was sufficient to consider the dimensions as follows.
Let $\eta=(p_1,p_2,\cdots, p_k)$ be a partition with $p_1 \geq p_2\geq\cdots\geq p_k$.
By $\hat{\eta}=(\hat{p}_1,\cdots,\hat{p}_m)$ denote the {\it dual partition} such that 
$\hat{p}_i\doteq \# \{j\vert p_j\geq i\}$. In the Young diagram language,
the diagram with {\it rows} consisting of $p_1,p_2,\cdots, p_k$ boxes, respectively
has {\it columns} consisting of $\hat{p}_1,\hat{p}_2,\cdots, \hat{p}_m$ boxes, respectively.
For example, the dual partition to $\eta=(5,3,3,1)$ is $\hat{\eta}=(4,3,3,1,1)$ as shows the Young diagram
of $\eta$
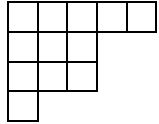
$$\yng(5,3,3,1)$$
\noindent Now, if $\eta=(p_1,p_2,\cdots, p_k)$ is a partition such that $p_1=t$ set 
\begin{equation}
n_1=\hat{p}_t; n_2= \hat{p}_{t-1}+\hat{p}_{t};\cdots; n_t= \hat{p}_1+\hat{p}_{2}+\cdots+\hat{p}_t.
\end{equation}
\noindent 
So $n_1,\cdots,n_t$ are the volumes of an increasing sequence of Young diagrams such that
the previous diagram is the result of collapsing the first column of the next one.
For example, the above partition yields the dimension vector $(1,2,5,8,12)$.
This way we define a vector $n(\eta)=(n_1,\cdots,n_t)$ and the set of
all such vectors can be characterized by the inequalities as follows:
\begin{equation}\label{kpineq}
n_1 \leq n_2 - n_1 \leq n_3 - n_2 \leq\cdots\leq n_t - n_{t-1}.
\end{equation}
\noindent In particular, this is a monotone sequence. Moreover, let $C$ be the Cartan matrix
of type $A_{t-1}$ and set $v=(n_1,\cdots,n_{t-1})$, $w=(0,\cdots,0,n_t)$. Then the formulae
(\ref{kpineq}) are equivalent to
\begin{equation}\label{nan0}
w - Cv\in {\bf Z}_+^{t-1}
\end{equation}
\begin{remark}\label{connect}
The condition (\ref{nan0}) has a very important sense in Nakajima's theory.
Namely, by \cite[Proposition~10.5]{na98} it is equivalent to the set ${\frak M}_0^{reg}(v,w)$ being
nonempty, which means that the generic orbit in $Z$ is closed with trivial stabilizer.
The most interesting general Nakajima's results hold under this condition and in this particular
case are just  equivalent to what is proved in \cite{kp}.
\end{remark}
A partition $\eta=(p_1,\cdots,p_k)$ of $t$ yields a nilpotent conjugacy class
$C_{\eta}$ of matrices with Jordan blocks of size $p_1,\cdots,p_k$, and
moreover, a special matrix $A\in C_{\eta}$ such that basis vectors
of ${\bf k}^t$ correspond to the boxes of Young diagram and $A$ maps
the boxes from the first column to 0 and each of the other boxes to its left neighbour. 

We now state a result from \cite{kp}, which is very close to Theorem \ref{7.2}.
Actually our statement is more strong than in \cite{kp} but one can easily check that
the original argument works for this statement without any change.
\begin{proposition}\label{quot} {\rm (Proposition 3.4 from \cite{kp})}

{\bf 1.} 
The map 
$\Theta:Z\to {\rm End}(U_t)$, $\Theta(A_1,B_1,\cdots,A_{t-1},B_{t-1})=A_{t-1}B_{t-1}$ 
is the categorical quotient with respect to $H$ for arbitrary dimension vector $(n_1,\cdots,n_t)$.

{\bf 2.} If $(n_1,\cdots,n_t)=n(\eta)$, then the image of $\Theta$ is equal $\overline{C_{\eta}}$.
\end{proposition}

\section{Nakajima's Theorem 7.3}

Before stating Nakajima's result we need some preliminary facts and notion.
Let $(n_1,\cdots,n_t)$ be a monotone dimension vector. Denote by ${\mathcal F}$
the variety of partial flags 
$\{0\}=E_0\subseteq E_1\subseteq E_2\subseteq\cdots\subseteq E_{t-1}\subseteq E_t={\bf k}^{n_t}$
with $\dim E_i=n_i$ for $i=1,\cdots,t$.
The variety ${\mathcal F}$ is projective and homogeneous with respect to the natural action
of $GL_{n_t}$, ${\mathcal F}\cong GL_{n_t}/P$, where $P$ is the stabilizer of a selected flag $f_0$,
a parabolic subgroup in $GL_{n_t}$. Recall that the tangent space
$T_{f_0} GL_{n_t}/P$ is isomorphic to ${\frak p}_0^*$, where ${\frak p}_0$ is the nilradical
of the Lie algebra of $P$.

Consider a closed subset $X\subseteq {\mathcal F}\times {\rm End}({\bf k}^{n_t})$ as follows:
\begin{equation}\label{def-x}
X=\{(f,A)\in {\mathcal F}\times {\rm End}({\bf k}^{n_t})\vert A E_i\subseteq E_{i-1},i=1,\cdots,t\}
\end{equation}
\noindent $X$ is naturally isomorphic to the cotangent bundle $T^*{\mathcal F}$
because the fiber of the projection $p_1:X\to {\mathcal F}$ over $f_0$
is $f_0\times {\frak p}_0$. 

Let $\mu$ be the dual partition to the ordered  sequence $(n_1,n_2-n_1,\cdots,n_t-n_{t-1})$ 
(in particular, if $(n_1,\cdots,n_t)=n(\eta)$, then $\mu =\eta$).
The following statement is well-known and can be found, e.g. in \cite[Theorem~3.3]{h}: 
\begin{proposition}\label{cruc}
$p_2(X)=\overline{C_{\mu}}$.
\end{proposition}

Now we need to introduce shortly quiver varieties in this particular case.
These are quotients by the action of a group, but
two papers, \cite{na94} and \cite{na98} propose two different approachs
to this notion, a K\"aler quotient and a quotient in the sense of Geometrical Invariant Theory,
respectively. Though the results we discuss are in \cite{na94}, we prefer the approach from \cite{na98}.

Nakajima considered two quotiens of $Z$ with respect to the action of $H$. The first,
${\frak M}_0$ is just the categorical quotient, ${\frak M}_0=Z/\!\!/H$ so the geometrical
points of ${\frak M}_0$ are in 1-to-1 correspondance with the closed $H$-orbits in $Z$.
On the other hand, one can
consider the semi-stable locus $Z^{ss}\subseteq Z$ (actually with respect to a particular choice of
a character of $H$ but we consider just one as in \cite{na98}). It is proved in \cite{na98} 
in general case that $Z^{ss}$ consists of stable points, that is, every $H$-orbit in $Z^{ss}$ 
is closed in $Z^{ss}$ and isomorphic to $H$. Hence,
there is a {\it geometric quotient} ${\frak M}=Z^{ss}/H$ (the construction of the quotient as
an algebraic variety is usual for GIT, see \cite[p.522]{na98}). In particular,
the points of ${\frak M}$ are in 1-to-1 correspondance with the $H$-orbits in $Z^{ss}$.
Moreover, the categorical quotient $Z\to {\frak M}_0$ gives rise to a natural map 
$\pi:{\frak M}\to{\frak M}_0$. Geometrically, $\pi$ sends a stable orbit $Hz$ to the unique 
closed (in $Z$) orbit in $\overline{Hz}$. Besides, the construction of ${\frak M}$ implies 
that $\pi$ is projective. Finally, Proposition \ref{quot} yields a convenient form
of $\pi$ as a map sending the stable orbit of $(A_1,B_1,\cdots,A_{t-1},B_{t-1})$ to
$A_{t-1}B_{t-1}\in{\rm End}(U_t)$.
\begin{theorem}{\rm (Theorem 7.3 from \cite{na94})} ${\frak M}\cong T^*{\mathcal F}$.
\end{theorem}
We want to reformulate the above theorem, as follows:
\begin{theorem}\label{my7.3}
There is an isomorphic map $\alpha:{\frak M}\to X\cong T^*{\mathcal F}$ making
the diagram commutative:

$$
\begin{CD}
{\frak M} @>\alpha>> X \\
@V{\pi}VV    @VVp_2V\\
\pi({\frak M}) @= \overline{C_{\mu}}
\end{CD}
$$
\end{theorem}

\begin{remark} 
Assume ${\frak M}_0^{reg}\neq\emptyset$, that is, the generic closed orbit in $Z$ is isomorphic to $H$.
By \cite[Proposition~3.24]{na98}, these generic closed orbits belong to $Z^{ss}$, hence,
$\pi({\frak M})$ is the whole of ${\frak M}_0$. 
We know that this happens precisely when the vector $(n_1,\cdots,n_t)$
is of Kraft-Procesi type (cf. Remark \ref{connect}). 
Reading the original proof by Nakajima, one can feel
that the author had the above diagram (with ${\frak M}_0$) in mind for
all monotone dimension vectors.
\end{remark}
We give a proof of Theorem \ref{my7.3} following the idea of the proof
of \cite[Theorem~7.3]{na94} but working in the setup of \cite{na98}, where we have:

\begin{proposition}\label{inject}
$(A_1,B_1,\cdots,A_{t-1},B_{t-1})\in Z^{ss}$ $\Leftrightarrow$
$A_1,\cdots,A_{t-1}$ are injective
\end{proposition}
\begin{proof}
In \cite[Lemma~3.8]{na98} we find a criterion of stability, which can be reformulated
in this case as follows: for any tuple of subspaces
$W_i\subseteq U_i,i=1,\cdots,t-1$ such that $A_i(W_i)\subseteq W_{i+1}$ and
$B_i(W_{i+1})\subseteq W_i$ for $i=1,\cdots,t-2$ and $A_{t-1}(W_{t-1})=0$
we have $W_1=W_2=\cdots W_{t-1}=0$.
Now assume that $A_1,\cdots,A_{p-1}$ are injective and $A_p$ is not.
Set $W_p={\rm Ker}(A_p)$ and $W_i=0$ for $i\neq p$. We claim that such a tuple contradicts
the above condition of stability. Indeed, we have $A_{p-1}B_{p-1}(W_p)=B_pA_p(W_p)=0$,
so $B_{p-1}(W_p)=0$, because $A_{p-1}$ is injective.
Conversely, assume that $A_1,\cdots,A_{t-1}$ are injective and let $W_1,\cdots,W_{t-1}$ be a
tuple as above. Since $A_{t-1}$ is injective and $A_{t-1}(W_{t-1})=0$, we have $W_{t-1}=0$.
Next, $A_{t-2}$ is injective and $A_{t-2}(W_{t-2})\subseteq W_{t-1}=0$ implies $W_{t-2}=0$.
So we get $W_1=W_2=\cdots W_{t-1}=0$ and the point is stable.
\end{proof}

\begin{proof} (of the Theorem)
Nakajima's construction for $\alpha$ is as follows:
\begin{equation}
\alpha((A_1,B_1,\cdots,A_{t-1},B_{t-1})=
\end{equation}
$$
=({\rm Im} A_{t-1}A_{t-2}\cdots A_1\subseteq 
{\rm Im} A_{t-1}A_{t-2}\cdots A_2\subseteq
\cdots \subseteq {\rm Im} A_{t-1} \subseteq U_t,
A_{t-1}B_{t-1})
$$
\noindent We claim that this map is well-defined. First of all
the components of $\alpha$ are $H$-invariant: this is clear for
the operator $A_{t-1}B_{t-1}$; as for the maps
$A_{t-1}\cdots A_i$ used to define the flag, the action
of $(h_1,\cdots,h_{t-1})\in H$ conjugates this map by the $h_i$
so does not change the image. Next, the constructed flag
belongs to ${\mathcal F}$ because by Proposition \ref{inject}
the maps $A_{t-1}\cdots A_i$ are injective over $Z^{ss}$ so the dimension
of the image is equal $n_i$. Finally, applying formulae (\ref{def_Z}) we have on $Z$: 
$$
A_{t-1}B_{t-1}A_{t-1}A_{t-2}\cdots A_i=
 A_{t-1}A_{t-2}B_{t-2}A_{t-2}\cdots A_i=
\cdots
$$
$$
\cdots= A_{t-1}\cdots A_{i+1}A_iB_iA_i
=A_{t-1}\cdots A_{i+1}A_{i}A_{i-1}B_{i-1}.
$$
\noindent Hence, the operator $A_{t-1}B_{t-1}$
maps ${\rm Im}A_{t-1}A_{t-2}\cdots A_i$ to  ${\rm Im}A_{t-1}A_{t-2}\cdots A_iA_{i-1}$.

Assume for $z=(A_1,B_1,\cdots, A_{t-1}, B_{t-1})$ and
$z'=(A'_1,B'_1,\cdots, A'_{t-1}, B'_{t-1})$: $z,z'\in Z^{ss}$ and $\alpha(z)=\alpha(z')$.
It is not difficult to see that we may conjugate $z'$ by an appropriate $h\in H$
such that not only the vector spaces ${\rm Im}A_{t-1}\cdots A_i$ and
${\rm Im}A'_{t-1}\cdots A'_i$ are equal for all $i$ but also $A_1=A'_1,A_2=A'_2,\cdots,
A_{t-1}=A'_{t-1}$. Then, applying the equality of the second parts of $\alpha$,
$A_{t-1}B_{t-1}=A'_{t-1}B'_{t-1}$ and having $A_{t-1}=A'_{t-1}$ is injective, 
we get $B_{t-1}=B'_{t-1}$. Next, we have 
\begin{equation}
A_{t-2}B_{t-2}=B_{t-1}A_{t-1}=B'_{t-1}A'_{t-1}=A'_{t-2}B'_{t-2}
\end{equation}
\noindent and $A_{t-2}=A'_{t-2}$ is injective, hence, $B_{t-2}=B'_{t-2}$.
Applying this argument repeatedly, we get $z=z'$, so $\alpha$ is injective.

On the other hand, for each point $x=((E_1\subseteq E_2\subseteq \cdots E_t),A)\in X$
we can identify $E_i$ with $U_i$, set $A_i$ to be the inclusion $E_i\subseteq E_{i+1}$
and set $B_i$ to be the restriction of $A$ to $E_{i+1}$. This way we get a point
$z\in Z$ and by Proposition \ref{inject} $z$ is stable. Since $\alpha(z)=x$, we proved that
$\alpha$ is bijective and moreover, $\alpha$ is an isomorphism, because $X$ is smooth.

The commutativity of the diagram follows from the definition of $\alpha$: indeed, we have $\pi=p_2\alpha$.
Finally, Proposition \ref{cruc} yields $p_2(X)=\overline{C_{\mu}}$.
\end{proof}

\begin{remark} A proof of \cite[Theorem~7.3]{na94} is also outlined in \cite{m}.
\end{remark}

Now we want to consider the action of $G= GL_{n_1}\times\cdots\times GL_{n_t}$ on $Z$.
Clearly, $Z^{ss}$ is $G$-stable and $G$ acts on both
$H$-quotients, ${\frak M}$ and ${\frak M}_0$ via the factor $GL_{n_t}$
such that $\pi$ is $G$-equivariant.
By Theorem \ref{my7.3} $\pi({\frak M})$ contains a dense $G$-orbit
so the same is true for ${\frak M}$. Since ${\frak M}$ is a geometric quotient $Z^{ss}/H$, we get:
\begin{corollary}
$Z^{ss}$ contains a dense $G$-orbit.
\end{corollary}
\begin{remark}
This corollary is exactly Theorem 2 from \cite{bhrr}. Indeed, 
in \cite[\S5]{bhrr} it is explained that the variety of representations of the Auslander
algebra in dimension ${\bf d}=(n_1,\cdots,n_t)$ is $Z$ (with the reverse numeration of the vector spaces).
Moreover, the $\Delta$-filtered
representations are all points in $Z$ with $A_1,\cdots,A_{t-1}$ being injective,
so by Proposition \ref{inject} this is $Z^{ss}$.
\end{remark}

In the next section we will see that $Z$ does not share these nice properties of $Z^{ss}$
\section{Nakajima's Theorem 7.2}\label{pf}
\begin{definition}\label{oper}
Let $\eta=(p_1,\cdots,p_s)$ be a partition of $n$ with $p_1\geq p_2\geq\cdots\geq p_s$.
For any $a\in {\bf Z}_+$ define a partition $\eta+a$ as follows:
If $a\geq s$, set $\eta+a=(p_1+1,\cdots,p_s+1,1,\cdots,1)$ so that
$\eta+a$ has $a-s$ more parts. Otherwise, if $s+a=2l$ set
$\eta+a=(p_1+1,\cdots,p_l+1,p_{l+1}-1,\cdots,p_s-1)$, or else,
if $s+a=2l+1$ set $\eta+a=(p_1+1,\cdots,p_l+1,p_{l+1},p_{l+2}-1,\cdots,p_s-1)$.
\end{definition}
Clearly, if $a\geq s$, then the Young diagam of $\eta+a$ is that of $\eta$ with added first
column of height $a$. 
Conversely, if $a < s$ the number of rows of $\eta+a$ can be less than that for $\eta$
(equal to $s$) provided $p_s=1$. As for the number of columns, it always increases by 1 from
$p_1$ to $p_1+1$. For example, if $\eta=(2,1,1)$, then
\begin{equation}\label{211}
\eta=\yng(2,1,1) \quad \eta+1=\yng(3,2) \quad \eta+2=\yng(3,2,1) \quad \eta+3=\yng(3,2,2)
\end{equation}

Recall that the set of partitions carries an  order  as follows:
\begin{equation}
\eta=(p_1,\cdots,p_s) \geq \nu=(q_1,\cdots,q_t) \Leftrightarrow
\sum_{i=1}^j p_i\geq \sum_{i=1}^j q_i,\;\forall j
\end{equation}
\noindent It follows from Definition \ref{oper}: 
\begin{equation}\label{transit}
\eta\geq\nu \Rightarrow \eta+a\geq\nu+a,\; \forall a\in{\bf Z}_+
\end{equation}
Consider the variety $L={\rm Hom}_{\bf k}({\bf k}^n,{\bf k}^{n+a})\times 
{\rm Hom}_{\bf k}({\bf k}^{n+a},{\bf k}^{n})$ of pairs $(A,B)$ of linear maps
and the maps $\pi:L\to {\rm End}({\bf k}^{n}), \pi(A,B)=BA$ and
$\rho:L\to {\rm End}({\bf k}^{n+a}), \rho(A,B)=AB$. The next lemma generalizes \cite[Lemma~2.3]{kp}:

\begin{lemma}\label{step}
$\rho(\pi^{-1}(\overline{C_{\eta}}))=\overline{C_{\eta+a}}$
\end{lemma}
\begin{proof}
To describe the pairs $(A,B)$ such that $\pi(A,B)=BA\in C_{\eta}$
we apply the techniques of so-called $ab$ diagrams from \cite[\S4]{kp}, as follows.
The $GL_n\times GL_{n+a}$-orbits of pairs $(A,B)$
such that $AB$ and $BA$ are nilpotent are depicted by the diagrams consisting of rows like $abab\cdots ab$,
which can start and end either with $a$ or with $b$. Clearly, the orbits are the isomorphism
classes of representations for the quiver with two vertices and two arrows of different directions,
so the diagram is nothing but a decomposition of the representation into indecomposable blocks
\begin{equation}\label{block}
e_1\xrightarrow{A}f_1\xrightarrow{B}e_2\xrightarrow{A}f_2\xrightarrow{B}\cdots
\xrightarrow{A}f_k\xrightarrow{B}0
\end{equation}
coresponding, for example, to the string $abab\cdots ab$ of length $2k$, where $e_1,\cdots,e_{k}$
and $f_1,\cdots,f_k$ are basis vectors of ${\bf k}^n$ and ${\bf k}^{n+a}$, respectively. Therefore,
the total number of letters $a$ in the diagram 
should be equal to  $n$ and that for $b$ is equal to $n+a$. We assume also that from the top to bottom
the rows are ordered by the length.
Now, it follows from (\ref{block}) that, taking all $a$-s from the diagram we get a partition of $n$, which corresponds to the conjugacy
class of $BA$ and the same for $b$-s and $AB$. Therefore we need to describe the $ab$ diagrams giving
$\eta$ as the $a$-part. To get such a diagram we have to fill in between of $a$-s the 
$p_1-1+p_2-1+\cdots+p_s-1$ letters $b$. After that we have $s+a$ more letters $b$ and we can add
them from the left and from the right to each row or create a row with single $b$.
In particular, if we want to place these $s+a$ letters $b$ as high as possible, then, 
if $a\geq s$ we add 2 $b$-s to each of $s$ rows and then add $a-s$ more rows with single $b$.
Otherwise, if $a < s$ and $s+a=2l$ we add 2 $b$-s to the first $l$ rows; else, if $s+a=2l+1$, 
we also add one $b$ to the $l+1$-th row. This way we get $\eta+a$ as the diagram of $AB$.
Of course, this is only one of possible ways to get an $ab$ diagram from $\eta$ but it is clear
from the above considerations that all other $b$-diagrams $\nu$ that we can get have the property
$\nu\leq \eta+a$. Moreover, we may take any other orbit $C_{\mu}\subseteq \overline{C_{\eta}}$
and the crucial property of the order is that $\mu\leq \eta$. So by the property (\ref{transit})
each $ab$-diagram over $\mu$ yields a $b$-diagram with partition less than $\mu+a$, hence,
than $\eta+a$ as well. So we proved that $\rho(\pi^{-1}(\overline{C_{\eta}}))$ contains $C_{\eta+a}$
and is contained in $\overline{C_{\eta+a}}$. On the other hand, by the First Fundamental Theorem
for $GL_n$ the map $\rho:L\to{\rm End}({\bf k}^{n+a})$ is the categorical quotient by $GL_n$.
Hence, the same is true for the restriction of $\rho$ to the closed $GL_n$-stable subvariety $\pi^{-1}(\overline{C_{\eta}})$. In particular, $\rho(\pi^{-1}(\overline{C_{\eta}}))$ is closed,
hence, is equal to $\overline{C_{\eta+a}}$.
\end{proof}

Now we are prepared to prove Theorem \ref{7.2} (Nakajima's Theorem 7.2):
\begin{proof}
By Proposition \ref{quot} we only need to show that the image
of $\Theta:Z\to {\rm End}(U_t), \Theta(A_1,B_1,\cdots,A_{t-1},B_{t-1})=A_{t-1}B_{t-1}$
is the closure of a nilpotent orbit and this follows from Lemma \ref{step}. 
Namely, we have $B_1A_1=0$ on $Z$, so we take $n=n_1, a = n_2-n_1$, and
$\eta=\rho=(1,1,\cdots,1)$. Hence, by Lemma \ref{step} the image of $A_1B_1$ is $\overline{C_{\eta+a}}$.
Then we take $n=n_1,a=n_3-n_2, \eta=\rho+(n_2-n_1)$, and get the image of
$A_2B_2$ to be $\overline{C_{(\rho+(n_2-n_1))+n_3-n_2}}$. Applying this argument repeatedly,
we complete the proof. 
\end{proof}
\begin{remark}
So the dense nilpotent conjugacy class in the image of $\Theta$ has the form of $C_{\lambda}$
with $\lambda=((\rho+(n_2-n_1))\cdots)+(n_t-n_{t-1})$. Only under Kraft-Procesi inequalities
on the dimension vector this partition has the clear direct connection with $(n_1,\cdots,n_t)=n(\lambda)$.
\end{remark}

As we already noted in the Introduction, the main part of results from \cite{kp} can not be generalized to the Nakajima's context. In particular, the following example shows that $Z$ can be reducible:
\begin{example}\label{reduc}
Take the dimension vector $(1,4,5)$. Then, applying the proof of Theorem \ref{7.2},
we have: $\rho=(1)$, the image of $A_1B_1$ is $\overline{C_{(2,1,1)}}$ and that
for $A_2B_2$ is $\overline{C_{(3,2)}}$, because $(1)+3=(2,1,1)$ and $(2,1,1)+1=(3,2)$ (c.f. (\ref{211})):
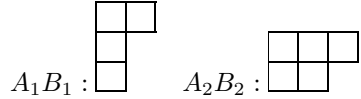
$$
A_1B_1: \yng(2,1,1) \quad A_2B_2: \yng(3,2)
$$ 
\noindent So as in the proof of Lemma \ref{step} we see that the pair $(A_2,B_2)$ coresponds to 
the $ab$-diagram $(babab,bab,a)$. 
Hence, $A_2$ is not injective (the summand $a$ corresponds
to the kernel of $A_2$). On the other hand, the stable representations in dimension $(1,4,5)$ exist
and, by Theorem \ref{my7.3} the image $\Theta(Z^{ss})=\overline{C_{(3,1,1)}}$.
Since the stability condition is open, if $Z$ would be irreducible, then $Z^{ss}$ would
be dense in $Z$, hence, $\Theta(Z^{ss})$ dense in $\Theta(Z)$, but we see this is false. 
\end{example}



\begin{thebibliography}{20}
\bibitem[BHRR]{bhrr} T. Br\"ustle, L. Hille, C.M. Ringel, and G. R\"ohrle,
The $\Delta$-filtered modules without self-extensions for the Auslander Algebra of $k[T]/\langle T^n\rangle$,
Algebras and Representation Theory {\bf 2} (1999), 295-312.
\bibitem[H]{h} W. Hesselink, Polarizations in the classical groups, Math. Z. {\bf 160} (1978), 217-234.
\bibitem[KP]{kp} H. Kraft and C. Procesi,
Closures of conjugacy classes of matrices are normal, Inventiones Math. 
{\bf 53} (1979), 227-247.
\bibitem[M]{m} A. Maffei, Quiver varieties of type {\it A}, Comment. Math. Helv. {\bf 80} (2005), 1, 1-27.
\bibitem[Na94]{na94} H. Nakajima, Instantons on ALE spaces, quiver varieties, 
and Kac-Moody algebras, Duke Math. J. {\bf 76} (1994) 2, 365–416. 
\bibitem[Na98]{na98} H. Nakajima, Quiver varieties and Kac-Moody algebras, Duke Math. J. {\bf 91} (1998),
3, 515-560.
\end{thebibliography}
\end{document}